\input amstex
\documentstyle{amsppt}
\magnification=\magstep1
 \hsize 13cm \vsize 18.35cm \pageno=1
\loadbold \loadmsam
    \loadmsbm
    \UseAMSsymbols
\topmatter
\NoRunningHeads
\title On the multiple $q$-Genocchi  and Euler numbers \endtitle
\author
  Taekyun Kim
\endauthor

\abstract The  purpose of this paper is to present a systemic
 study of some families of multiple $q$-Genocchi and Euler numbers by using multivariate
 $q$-Volkenborn integral (= $p$-adic $q$-integral) on $\Bbb
 Z_p $. From the studies of those $q$-Genocchi numbers and
 polynomials of higher order we derive some interesting identities related
 to $q$-Genocchi numbers and polynomials of higher order.

\endabstract
\thanks  2000 AMS Subject Classification: 11B68, 11S80
\newline  This paper is supported by  Jangjeon Research Institute for Mathematical
Science(JRIMS-10R-2001)
\endthanks
\endtopmatter

\document

{\bf\centerline {\S 1. Introduction}}

 \vskip 20pt

Let $p$ be a fixed odd prime. Throughout this paper $\Bbb Z_p ,$
$\Bbb Q_p ,$ $\Bbb C,$ and $\Bbb C_p$ will, respectively, denote the
ring of $p$-adic rational integers, the field of $p$-adic rational
numbers, the complex number field, and the completion of algebraic
closure of $\Bbb Q_p .$ Let $v_p$ be the normalized exponential
valuation of $\Bbb C_p$ with $|p|_p=p^{-v_p(p)}=p^{-1}$ and let $q$
be regarded as either a complex number $q\in\Bbb C$ or a $p$-adic
number $q\in\Bbb C_p$. If $q\in\Bbb C$, then we always assume
$|q|<1$. If $q\in\Bbb C_p$, we normally assume
$|1-q|_p<p^{-\frac{1}{p-1}},$ which implies that $q^x =\exp(x\log q
) $ for $|x|_p\leq 1$. Here, $| \cdot |_p$ is the $p$-adic absolute
value in $\Bbb C_p$ with $|p|_p=\frac{1}{p}.$ The $q$-basic natural
number are defined by $[n]_q=\frac{1-q^n}{1-q}=1+q+\cdots+q^{n-1}, $
( $n\in \Bbb N $), and $q$-factorial are also defined  as $[n]_q
!=[n]_q \cdot [n-1]_q \cdots [2]_q \cdot [1]_q.$ In this paper we
use the notation of Gaussian binomial coefficient as follows:
$${\binom {n}{k}}_q =\frac{[n]_q!}{[n-k]_q![k]_q!}
=\frac{[n]_q\cdot [n-1]_q \cdots [n-k+1]_q}{[k]_q!}. \tag1$$ Note
that $\lim_{q\rightarrow 1}{\binom{n}{k}}_q=\binom{n}{k}=\frac{n
\cdot(n-1)\cdots(n-k+1)}{n!}.$ The Gaussian coefficient satisfies
the following recursion formula:
$$\binom{n+1}{k}_q =\binom{n}{k-1}_q+q^k\binom{n}{k}_q=q^{n-k}\binom{n}{k-1}_q+\binom{n}{k}_q, \text{ cf. [12].} \tag2$$
From thus recursion formula we derive
$$\binom{n}{k}_q=\sum_{d_0+\cdots+d_k=n-k, d_i\in\Bbb N} q^{d_1+2d_2+\cdots+kd_k},\text{ cf.[1, 2, 12, 13, 14].}$$
The $q$-binomial formulae are known as
$$(b;q)_n=\prod_{i=1}^n\left(1-bq^{i-1}\right)=\sum_{k=0}^n\binom{n}{k}_qq^{\binom{k}{2}}(-1)^kb^k
,$$ and
$$\frac{1}{(b;q)_n}=\prod_{i=1}^n\left(1-bq^{i-1}\right)^{-1}=\sum_{k=0}^{\infty}\binom{n+k-1}{k}_q
b^k, \text{ cf.[12]}. \tag3 $$ In this paper we use the notation
$$[x]_q=\frac{1-q^x}{1-q}, \text{ and  }
[x]_{-q}=\frac{1-(-q)^x}{1+q}.$$ Hence, $\lim_{q\rightarrow
1}[x]_q =1, $ for any $x$ with $|x|_p\leq 1$ in the present
$p$-adic case, cf.[1-18].

 For $d$ a fixed positive
integer with $(p,d)=1$, let
$$\split
& X=X_d = \lim_{\overleftarrow{N} } \Bbb Z/ dp^N \Bbb Z , \ X_1 =
\Bbb Z_p , \cr  & X^\ast = \underset {{0<a<d p}\atop {(a,p)=1}}\to
{\cup} (a+ dp \Bbb Z_p ), \cr & a+d p^N \Bbb Z_p =\{ x\in X | x
\equiv a \pmod{dp^n}\},\endsplit$$ where $a\in \Bbb Z$ lies in
$0\leq a < d p^N$. In [9], we note  that
$$\mu_{-q}(a+dp^N \Bbb
Z_p)=(1+q)\frac{(-1)^aq^a}{1+q^{dp^N}}
=\frac{(-q)^a}{[dp^N]_{-q}},$$ is distribution on $X$ for $q\in \Bbb
C_p$ with $|1-q|_p< 1. $ We say that $f$ is a uniformly
differentiable function at a point $a \in\Bbb Z_p $ and denote this
property by $f\in UD(\Bbb Z_p )$, if the difference quotients $F_f
(x,y) = \dfrac{f(x) -f(y)}{x-y} $ have a limit $l=f^\prime (a)$ as
$(x,y) \to (a,a)$. For $f\in UD(\Bbb Z_p )$, this distribution
yields an integral as follows:
$$I_{-q}=\int_{\Bbb Z_p} f(x)d\mu_{-q}(x)=\int_{X} f(x)d\mu_{-q}(x)=\lim_{N\rightarrow
\infty}\frac{1}{[dp^N]_{-q}}\sum_{x=0}^{dp^N-1}f(x)(-q)^x,
$$
which has a sense as we see readily that the limit is convergent
(see [9, 10, 14, 15]). Let $q=1$. Then we have the fermionic
$p$-adic integral on $\Bbb Z_p$ as follows:
$$I_{-1}=\int_{\Bbb
Z_p}f(x)d\mu_{-1}(x)=\lim_{N\rightarrow
\infty}\sum_{x=0}^{p^N-1}f(x)(-1)^x, \text{ cf.[3, 6, 7, 8,  9,
13, 14].}$$
 For any positive integer $N,$ we set
$$\mu_q (a+lp^N \Bbb Z_p)=\frac{q^a }{[lp^N ]_q}, \text{ cf. [5, 9, 15, 16, 17, 18]},$$
and this can be extended to a distribution on $X$. This
distribution yields $p$-adic bosonic $q$-integral as follows (see
[12, 17, 18]):
$$ I_q (f)=\int_{\Bbb Z_p}  f(x) \,d\mu_q(x)=\int_X  f(x) \,d \mu_q (x) ,$$
where $f \in UD(\Bbb Z_p )=\text{ the space of uniformly
differentiable function  on } \Bbb Z_p $ with values in $\Bbb
C_p$.

In view of notation, $I_{-1}$ can be written symbolically as
$I_{-1}(f)=\lim_{q\rightarrow -1} I_{q}(f),$ cf.[9].  For
$n\in\Bbb N$, let $f_n(x)=f(x+n)$. Then we have
$$q^nI_{-q}(f_n)=(-1)^nI_{-q}(f)+[2]_{q}\sum_{l=0}^{n-1}(-1)^{n-1-l}q^l f(l), \text{ see [9]}.\tag4$$
For any complex number $z$, it is well known that the familiar Euler
polynomials $E_n(z)$ are defined by means of the following
generating function:
$$F(z,t)=\frac{2}{e^{t}+1}e^{zt}=\sum_{n=0}^{\infty}E_n(z)\frac{t^n}{n!},
\text{ for $|t|<\pi $, cf.[13,14]}.$$ We note that, by
substituting $z=0$, $E_n(0)=E_n$ is the familiar $n$-th Euler
number defined by
$$F(t)=F(0,t)=\frac{2}{e^t+1}=\sum_{n=0}^{\infty}E_n\frac{t^n}{n!}, \text{ cf.[12]}.$$
The Genocchi  numbers $G_n$ are defined by the generating function
$$\frac{2t}{e^t+1}=\sum_{n=0}^{\infty}G_n\frac{t^n}{n!},
(|t|<\pi).$$ It satisfies $G_1=1, G_3=G_5=\cdots=G_{2k+1}=0,$ and
even coefficients are given by
$$G_n=2(1-2^n)B_n=2nE_{2n-1}(0),$$
where $B_n$ are Bernoulli numbers and $E_n(x)$ are Euler
polynomials. By meaning of the generalization of $E_n$,
Frobenius-Euler numbers and polynomials are also defined by
$$\frac{1-u}{e^t-u}=\sum_{n=0}^{\infty}H_n(u)\frac{t^n}{n!}, \text{
and }
\frac{1-u}{e^t-u}e^{xt}=\sum_{n=0}^{\infty}H_n(u,x)\frac{t^n}{n!},
\text{ for $u\in\Bbb C$, cf.[12, 14]}.$$ Over five decades ago,
Carlitz [1, 2] defined $q$-extension of Frobenius-Euler numbers
and polynomials and proved properties analogous to those satisfied
$H_n(u) $ and $H_n(u,x)$. In previous my paper [6, 7, 8] the
author defined the q-extension of ordinary Euler and polynomials
and proved properties analogous to those satisfied $E_n$ and
$E_n(x)$. In [6] author also constructed the $q$-Euler numbers and
polynomials of higher order and gave some interesting formulae
related to Euler numbers and polynomials of higher order . The
purpose of this paper is to present a systemic
 study of some families of multiple $q$-Genocchi and Euler numbers by using multivariate
 $q$-Volkenborn integral (= $p$-adic $q$-integral) on $\Bbb
 Z_p $. From the studies of these $q$-Genocchi numbers and
 polynomials of higher order we derive some interesting identities related to
  $q$-Genocchi numbers and polynomials of higher order.

\vskip 20pt

{\bf\centerline {\S 2. Preliminaries / $q$-Euler polynomials}}
\vskip 10pt In this section we assume that $q\in\Bbb C_p$ with
$|1-q|_p<p^{-\frac{1}{p-1}}$. Let $f_1(x)$ be translation with
$f_1(x)=f(x+1)$. From (4) we can derive
$$ I_{-1}(f_{1})=I_{-1}(f)+2f(0).$$
If we take $f(x)=e^{(x+y)t}$, then we have Euler polynomials from
the integral equation of $I_{-1}(f)$ as follows:
$$  \int_{\Bbb Z_p}e^{(x+y)t}d\mu_{-1}(y)=e^{xt}\frac{2}{e^{t}+1}=\sum_{n=0}^{\infty}\frac{E_n(x)t^n}{n!}.$$
That is,
$$\int_{\Bbb Z_p}y^nd\mu_{-1}(y)=E_n, \text{ and} \int_{\Bbb Z_p}(x+y)^nd\mu_{-1}(y)=E_{n}(x).$$
Now we consider the following  multivariate $p$-adic fermionic
integral on $\Bbb Z_p$ as follows:
$$\int_{\Bbb Z_p}\cdots \int_{\Bbb Z_p} e^{(x_1 +\cdots +x_r+x)t}d\mu_{-1}(x_1) \cdots d\mu_{-1}(x_r)
=\left(\frac{ 2}{e^t+1}\right)^r e^{xt}=\sum_{n=0}^{\infty}
E_n^{(r)}(x)\frac{t^n}{n!},\tag5$$ where $E_n^r(x)$ are the Euler
polynomials of order $r$.

From (5) we note that
$$\int_{\Bbb Z_p}\cdots \int_{\Bbb Z_p} (x_1 +\cdots +x_r+x)^nd\mu_{-1}(x_1) \cdots
d\mu_{-1}(x_r)=E_n^{(r)}(x). \tag6$$ In view of (6) we can define
the q-extension of Euler polynomials of higher order.  For $h\in\Bbb
Z, $ $k\in\Bbb N$, let us consider the extended higher order
$q$-Euler polynomials as follows:
$$E_{m,q}^{(h,k)}(x)=\int_{\Bbb Z_p} \cdots \int_{\Bbb Z_p}[x_1+ \cdots
+ x_k +x]_q^m q^{\sum_{j=1}^k x_j (h-j)}d\mu_{-q}(x_1)\cdots
d\mu_{-q}(x_r), \text{ see [6] }.$$ From this definition we can
derive
$$E_{m,q}^{(h,k)}(x)=[2]_q^k\frac{1}{(1-q)^m}\sum_{j=0}^m\binom{m}{j}(-1)^jq^{xj}\frac{1}{(-q^{j+h};q^{-1})_k}
, \text{ see [6] }.\tag7$$ It is easy to see that
$$(-q^{j+k-1};q^{-1})_k=\prod_{l=0}^{k-1}(1+q^lq^j)=(-q^j;q)_k $$
In the special case $h=k-1$, we can easily see that
$$\aligned
E_{m,q}^{(k-1.k)}(x)&=[2]_q^k\frac{1}{(1-q)^m}\sum_{j=0}^m\binom{m}{j}(-1)^jq^{xj}\frac{1}{(-q^{j+k-1};q^{-1})_k}\\
&=[2]_q^k\frac{1}{(1-q)^m}\sum_{j=0}^m\binom{m}{j}(-1)^jq^{xj}\frac{1}{(-q^j;q)_k}\\
&=[2]_q^k\frac{1}{(1-q)^m}\sum_{j=0}^m\binom{m}{j}(-1)^jq^{xj}\sum_{n=0}^{\infty}{\binom{k+n-1}{n}}_q(-q^j)^n\\
&=[2]_q^k\sum_{n=0}^{\infty}{\binom{k+n-1}{n}}_q(-1)^n[n+x]_q^m.
\endaligned\tag8$$
Let $F^k(t,x)=\sum_{n=0}^{\infty}E_{n,q}^{(k-1,k)}(x)$ be the
generating function of $E_{n,q}^{k-1,k}(x)$. From (8) we note that
$$\aligned
F_q^k(t,x)&=\sum_{m=0}^{\infty}E_{m,q}^{(k-1,k)}(x)\frac{t^m}{m!}=[2]_q^k\sum_{m=0}^{\infty}
\sum_{n=0}^{\infty}{\binom {k+n-1}{n}}_q(-1)^n[n+x]_q^m\frac{t^m}{m!}\\
&=[2]_q^k\sum_{n=0}^{\infty}{\binom {k+n-1}{n}}_q(-1)^n e^{[n+x]_q}.
\endaligned$$

Remark. For $ w\in \Bbb C_p$ with $|1-w|<1$, we have
$$I_{-1}(w^xe^{tx})=\int_{\Bbb Z_p}w^xe^{tx}d\mu_{-1}(x)=\frac{2}{w
e^t+1}=\sum_{n=0}^{\infty}E_n(w)\frac{t^n}{n!}, \text{ see [3, 9,
11, 12, 15] }.$$ Thus, we note that $\int_{\Bbb Z_p}w^x x^n
d\mu_{-1}(x)=E_{n}(w) $ and $E_n(w)=\frac{2}{w+1}H_n(-w^{-1}), $
where $H_n(-w^{-1})$ are Frobenius-Euler numbers.

In the previous paper [11], the $q$-extension of $E_{n}(w)$(=
twisted $q$-Euler numbers) are studied as follows:
$$I_{-q}(w^xe^{t[x]_q})=\int_{\Bbb
Z_p}w^xe^{t[x]}d\mu_{-q}(x)=\sum_{n=0}^{\infty}E_{n,q}(w)\frac{t^n}{n!}
. \tag9$$ From (9) we note that
$$E_{n,q}(w)=\int_{\Bbb
Z_p}w^x[x]_q^n
d\mu_{-q}(x)=\frac{[2]_q}{(1-q)^n}\sum_{j=0}^{n}\binom{n}{j}(-1)^j\frac{1}{1+q^{j+1}w},
\text{ see [11] }. $$ By the exactly same method of Eq.(7), we can
also derive the multiple twisted $q$-Euler numbers as follows:
$$\aligned
&E_{m,q}^{(h,k)}(w,x)\\
&=\int_{\Bbb Z_p}\cdots\int_{\Bbb Z_p}w^{\sum_{j=1}^k x_j}[x_1 +
\cdots + x_k+x]_q^mq^{\sum_{j=1}^k (h-j)x_j}  d\mu_{-q}(x_1)\cdots
d\mu_{-q}(x_k) .\endaligned \tag10 $$ From (10) we can easily derive
$$E_{m,q}^{(h,k)}(w,x)=\frac{[2]_q^k}{(1-q)^m}\sum_{l=0}^m\binom{m}{l}(-q^x)^l\frac{1}{(-wq^{h+l};q^{-1})}.$$
For $h=k-1$, we have
 $$\aligned
  E_{m,q}^{(k-1,k)}(w,x)&=\frac{[2]_q^k}{(1-q)^m}\sum_{l=0}^m\binom{m}{l}(-q^x)^l\frac{1}{(-wq^{l};q)_k}\\
& =[2]_q^k\sum_{n=0}^{\infty}{\binom{n+k-1}{n}}_q(-w)^n[n+x]_q^m.
 \endaligned \tag11$$
Let $F_q^k(w,x)=\sum_{m=0}^{\infty}E_{m,q}^{(k-1,
k)}(x)\frac{t^m}{m!}$. From (11), we can easily derive
$$F_q^k(w,x)=[2]_q^k\sum_{n=0}^{\infty}{\binom{n+k-1}{n}}_q(-w)^ne^{[n+x]_qt}.
$$
Remark. When we consider those $q$-Euler numbers and polynomials in
complex number field, we assume that $q\in\Bbb C$ with $|q|<1$.

\vskip 10pt

 {\bf\centerline {\S 3. Genocchi and $q$-Genocchi numbers   }}
  \vskip 10pt

From (4) we note that
$$t\int_{\Bbb Z_p}e^{xt}d\mu_{-1}(x)=\frac{2t}{e^t+1}=\sum_{n=0}^{\infty}G_n\frac{t^n}{n!}. \tag12$$
Thus, we have
$$\int_{\Bbb
Z_p}e^{xt}d\mu_{-1}(x)=\sum_{n=0}^{\infty}\frac{G_{n+1}}{n+1}\frac{t^n}{n!}.\tag13$$
By (13) we easily see that
$$\int_{\Bbb Z_p}x^nd\mu_{-1}(x)=\frac{G_{n+1}}{n+1}, \text{ and }
\int_{\Bbb Z_p}
(x+y)^nd\mu_{-1}d\mu_{-1}(y)=\frac{G_{n+1}(x)}{n+1},$$ where
$G_n(x)$ are Genocchi polynomials ( see [8] ).

From the multivariate $p$-adic fermionic integral on $\Bbb Z_p$ we
can also derive the Genocchi numbers of order $r$ as follows:
$$t^r\int_{\Bbb Z_p}\cdots \int_{\Bbb
Z_p}e^{(x_1 + \cdots + x_r )t}d\mu_{-1}(x_1)\cdots d\mu_{-1}(x_{r})
=\left(\frac{2t}{e^t+1}\right)^r=\sum_{n=0}^{\infty}G_n^{(r)}\frac{t^n}{n!},
\text{ $r\in\Bbb N,$} \tag 14$$ where $G_n^{(r)}$ are the Genocchi
numbers of order $r$.

From (14) we note that
$$\sum_{n=0}^{\infty}
\int_{\Bbb Z_p}\cdots \int_{\Bbb Z_p}(x_1+\cdots+x_r)^n
d\mu_{-1}(x_1)\cdots d\mu_{-1}(x_{r})\frac{(r+n)_rt^{n+r}}{(n+r)!}
=\sum_{n=0}^{\infty}G_n^{(r)}\frac{t^n}{n!},\tag15$$ where
$(x)_r=x(x-1)\cdots (x-r+1) .$ By (14) and (15), we easily see that
$$\int_{\Bbb Z_p}\cdots \int_{\Bbb Z_p}(x_1+\cdots+x_r)^n
d\mu_{-1}(x_1)\cdots
d\mu_{-1}(x_{r})=\frac{1}{r!\binom{n+r}{r}}G_{n+r}^{(r)}, \text{
where $n\in\Bbb N\cup \{0\}$},\tag16 $$ and
$G_0^{(r)}=G_1^{(r)}=\cdots=G_{r-1}^{(r)}=0.$ Thus, we obtain the
following theorem.

\proclaim{ Theorem 1} For $n\in\Bbb N\cup \{0\}$, $r\in\Bbb N$, we
have
$$\aligned
G_{n+r}^{(r)}&=(n+r)_r\int_{\Bbb Z_p}\cdots \int_{\Bbb Z_p}(x_1
+\cdots +x_r)^n d\mu_{-1}(x_1)\cdots d\mu_{-1}(x_r)\\
&=\binom{n+r}{r}r!\int_{\Bbb Z_p}\cdots \int_{\Bbb Z_p}(x_1 +\cdots
+x_r)^n d\mu_{-1}(x_1)\cdots d\mu_{-1}(x_r),
\endaligned$$
where $(x)_r=x(x-1)\cdots (x-r+1) .$
\endproclaim
Recently we constructed the $q$-extension of Genocchi numbers as
follows:
$$t\int_{\Bbb
Z_p}e^{[x]_qt}d\mu_{-q}(x)=[2]_qt\sum_{n=0}^{\infty}(-1)^nq^n
e^{[n]_qt}=\sum_{n=0}^{\infty}G_{n,q}\frac{t^n}{n!}, \text{ see
[8] }. \tag17$$

Thus, we note that
$$\int_{\Bbb Z_p}[x]_q^n
d\mu_{-q}(x)=G_{n,q}=n\frac{[2]_q}{(1-q)^{n-1}}\sum_{l=0}^{n-1}\binom{n-1}{l}\frac{(-1)^l}{1+q^{l+1}},
\text{ see [8] }.$$

In view of (14) we can define the q-extension of Genocchi numbers of
higher order. For $h\in\Bbb Z, $ $k\in\Bbb N$, let us consider the
extended higher order $q$-Genocchi numbers  as follows:

$$\aligned
&\sum_{n=0}^{\infty}G_{n,q}^{(h,k)}\frac{t^n}{n!}=t^k\int_{\Bbb
Z_p}\cdots\int_{\Bbb Z_p}e^{[x_1+\cdots+x_k]_qt}q^{\sum_{j=1}^k
x_{j}(h-j)}d\mu_{-q}(x_{1}) \cdots d\mu_{-q}(x_{k})\\
&=\sum_{n=0}^{\infty}\int_{\Bbb Z_p}\cdots\int_{\Bbb
Z_p}[x_1+\cdots+x_k]_q^nq^{\sum_{j=1}^k x_{j}(h-j)}d\mu_{-q}(x_{1})
\cdots d\mu_{-q}(x_{k})\frac{t^{n+k}}{n!}\\
&=\sum_{n=0}^{\infty}\int_{\Bbb Z_p}\cdots\int_{\Bbb
Z_p}[x_1+\cdots+x_k]_q^nq^{\sum_{j=1}^k x_{j}(h-j)}d\mu_{-q}(x_{1})
\cdots d\mu_{-q}(x_{k})\frac{(n+k)_k t^{n+k}}{(n+k)!}.
\endaligned$$
Thus, we have
$$\aligned
G_{n+k,q}^{(h,k)}&=k!\binom{n+k}{k}\int_{\Bbb Z_p}\cdots\int_{\Bbb
Z_p}[x_1+\cdots+x_k]_q^nq^{\sum_{j=1}^k x_{j}(h-j)}d\mu_{-q}(x_{1})
\cdots d\mu_{-q}(x_{k})\\
&=k!\binom{n+k}{k}\frac{[2]_q^k}{(1-q)^n}\sum_{l=0}^{n}\binom{n}{l}(-1)^l\frac{1}{(-q^{h+l};q^{-1})_k}
, \endaligned$$ and
$$G_{0,q}^{(h,k)}=G_{1,q}^{(h,k)}=\cdots=G_{k-1,q}^{(h,k)}=0.$$
If we take $h=k-1$, then we have
$$\aligned
G_{n+k,q}^{(k-1,k)}&=k!\binom{n+k}{k}\frac{[2]_q^k}{(1-q)^n}\sum_{l=0}^n\binom{n}{l}(-1)^l\frac{1}{(-q^l;q)_k}\\
&=k!\binom{n+k}{k}\frac{[2]_q^k}{(1-q)^n}\sum_{l=0}^n\binom{n}{l}(-1)^l\sum_{m=0}^{\infty}{\binom{m+k-1}{m}}_q(-q^l)^m\\
&=k!\binom{n+k}{k}[2]_q^k\sum_{m=0}^{\infty}{\binom{m+k-1}{m}}_q(-1)^m[m]_q^n.
\endaligned$$
Therefore we obtain the following theorem.

\proclaim{ Theorem 2} For $h\in\Bbb Z$, $k\in\Bbb N$, we have
$$G_{n+k,q}^{(h,k)}=k!\binom{n+k}{k}\frac{[2]_q^k}{(1-q)^n}\sum_{l=0}^n\binom{n}{l}(-1)^l\frac{1}{(-q^{h+l};q^{-1})_k},$$
and
$$G_{n+k,q}^{(k-1,k)}=k!\binom{n+k}{k}[2]_q^k\sum_{m=0}^{\infty}{\binom{m+k-1}{m}}_q(-1)^m[m]_q^n.$$
\endproclaim

Let
$$h_q^k(t)=\sum_{n=0}^{\infty}G_{n,q}^{(k-1,k)}\frac{t^n}{n!}
=\sum_{n=0}^{\infty}G_{n+k,q}^{(k-1,k)}\frac{t^{n+k}}{(n+k)!},\tag18$$
because $G_{0,q}^{(k-1,k)}=\cdots=G_{k-1,q}^{(k-1,k)}=0 .$ By (18)
and Theorem 2, we see that
$$\aligned
h_q^k(t)&=\sum_{n=0}^{\infty}G_{n,q}^{(k-1,k)}\frac{t^n}{n!}
=[2]_q^kt^k\sum_{n=0}^{\infty}\sum_{m=0}^{\infty}{\binom{m+k-1}{m}}_q(-1)^m[m]_q^n\frac{t^m}{m!}\\
&=[2]_q^kt^k\sum_{m=0}^{\infty}{\binom{m+k-1}{m}}_q(-1)^me^{[m]_qt}.
\endaligned$$

Remark. For $w\in\Bbb C_p$ with $|1-w|_p<1$, we can also define
$w$-Genocchi numbers (= twisted Genocchi numbers) as follows:
 $$t\int_{\Bbb
 Z_p}w^xe^{xt}d\mu_{-1}(x)=\frac{2t}{we^t+1}=\sum_{n=0}^{\infty}G_{n,w}
 \frac{t^n}{n!}, \text{ cf.[3, 11, 13] }. $$
From this we note that $\lim_{w\rightarrow 1}G_{n,w}=G_n.$ The
$q$-extension of $G_{n,w}$ can be also defined  by
$$t\int_{\Bbb Z_p}w^x e^{[x]_q t}
d\mu_{-q}(x)=\sum_{n=0}^{\infty}G_{n,q,w}\frac{t^n}{n!}, \text{
cf. [3, 8, 11, 13] }. \tag19$$ By (19) we easily see that
$$G_{n,q,w}=n\frac{[2]_q}{(1-q)^{n-1}}
\sum_{l=0}^{n-1}\binom{n-1}{l}\frac{(-1)^l}{1+q^{l+1}w}, \text{
cf.[11] }.$$ From this we also note that $\lim_{w\rightarrow
1}G_{n,q,w}=G_{n,q}. $

Now we consider the extended $(q,w)$-Genocchi numbers by using
multivariate $p$-adic  fermionic integral on $\Bbb Z_p$. For
$h\in\Bbb Z, k\in \Bbb N,$ $w\in \Bbb C_p$ with $|1-w|_p<1,$ we
define $G_{n,q, w}^{(h,k)}$ as follows:
$$\aligned
&\sum_{n=0}^{\infty}G_{n,q,w}^{(h,k)}\frac{t^n}{n!} \\
&=t^k\int_{\Bbb Z_p}\cdots \int_{\Bbb Z_p}w^{x_1 +\cdots +x_k }
e^{[x_1+\cdots+x_k]_qt}q^{\sum_{j=1}^k
x_j(h-j)}d\mu_{-q}(x_{1})\cdots d\mu_{-q}(x_{k}).\endaligned
\tag20$$ From (20) we can derive
$$G_{n+k,q,w}^{(h,k)}=k!\binom{n+k}{k}\frac{[2]_q^k}{(1-q)^n}\sum_{l=0}^n
\binom{n}{l}(-1)^l\frac{1}{(-wq^{h+l};q^{-1})_k}, $$ and
$$G_{n+k,q,w}^{(k-1,k)}=k!\binom{n+k}{k}[2]_q^k
\sum_{m=0}^{\infty}{\binom{m+k-1}{m}}_q(-w)^m [n]_q^m.$$ Let
$h_{q,w}^k(t)=\sum_{n=0}^{\infty}G_{n,q}^{(k-1,k)}\frac{t^n}{n!}. $
Then we have
$$ h_{q,w}^k(t)=\sum_{n=0}^{\infty}G_{n+k,q}^{(k-1,
k)}\frac{t^{n+k}}{(n+k)!}=[2]_q^k
t^k\sum_{m=0}^{\infty}{\binom{m+k-1}{m}}_q (-w)^m e^{[m]_qt}. $$

\vskip 20pt

\quad Taekyun Kim

\quad EECS, Kyungpook National University, Taegu 702-701, S. Korea

\quad e-mail:\text{ tkim$\@$knu.ac.kr; tkim64$\@$hanmail.net}

\enddocument